\documentclass[titlepage,11pt]{article}
\usepackage{latexsym}
\usepackage{amsfonts}
\usepackage{amsmath}
\usepackage{hyperref}
\hypersetup{
    pdftitle=   {A relative of Hadwiger's conjecture},
   pdfauthor=  {Katherine Edwards, Dong Yeap Kang, Jaehoon Kim, Sang-il Oum, Paul Seymour}
}

\newcommand\blackslug{\hbox{\hskip 1pt \vrule width 4pt height 8pt depth 1.5pt
        \hskip 1pt}}
\newcommand\bbox{\hfill \quad \blackslug \bigbreak}

\def\l{,\ldots,}

\title{A relative of Hadwiger's conjecture}
\author{Katherine Edwards\thanks{Supported by an NSERC PGS-D3 Fellowship and a Gordon Wu Fellowship.}, 
Princeton University, Princeton, NJ 08544, USA
\\
\\
Dong Yeap Kang\footnotemark[2], KAIST, Daejeon, 34141, Republic of Korea
\\
\\
Jaehoon Kim\footnotemark[2], University of Birmingham, Birmingham, B15 2TT, UK
\\
\\
Sang-il Oum\thanks{Supported by Basic Science Research
  Program through the National Research Foundation of Korea (NRF)
  funded by  the Ministry of Science, ICT \& Future Planning
  (2011-0011653).
}, KAIST, Daejeon, 34141, Republic of Korea
\\
\\
Paul Seymour\thanks{Supported by ONR grant N00014-10-1-0680 and 
NSF grant DMS-1265563.}, Princeton University, Princeton, NJ 08544, USA}

\date{July 18, 2014; revised August 21, 2015}

\newtheorem{thm}{}[section]

\newcommand{\Proof}{\noindent{\bf Proof.}\ \ }

\begin{document}
\maketitle
\begin{abstract}
Hadwiger's conjecture asserts that if a simple graph $G$ has no $K_{t+1}$ minor,
then its vertex set $V(G)$ can be partitioned into $t$ stable sets. This is 
still open, but we prove under the same hypotheses that $V(G)$ can be partitioned into 
$t$ sets $X_1\l X_t$, such that for $1\le i\le t$, the subgraph induced on $X_i$ has maximum degree at
most a function of $t$. This is sharp, in that the conclusion
becomes false if we ask for a partition into $t-1$ sets with the same property.
\end{abstract}

\section{Introduction}
All graphs in this paper are finite and have no loops or multiple edges. A graph $H$ is a {\em minor}
of a graph $G$ if a graph isomorphic to $H$ can be obtained from a subgraph of $G$ by edge-contraction.
In 1943, Hadwiger~\cite{hadwiger} proposed the following, perhaps the most famous open problem in graph theory:
\begin{thm}\label{hadwiger}
{\bf (Hadwiger's Conjecture.)} For all integers $t\ge 0$, and every graph $G$, if $K_{t+1}$ is not a minor of $G$,
then the chromatic number of $G$ is at most $t$; that is, $V(G)$ can be partitioned%
\footnote{A  \emph{partition} of a set $V$ is a list of pairwise disjoint (possibly empty) subsets of $V$ whose union equals $V$.}
 into $t$ stable sets.
\end{thm}
This remains open, although it has been proved for all $t\le 5$ (see~\cite{RST}). 
It is best possible in that the result becomes false
if we ask for a partition into $t-1$ stable sets.

There are several results proving weakenings of Hadwiger's conjecture (see section~3), but as far as we know, the
result of this paper is the first which (under the same hypotheses as \ref{hadwiger}) asserts the existence of a partition of $V(G)$ 
into $t$ sets with any non-trivial property. We prove
the following. (If $G$ is a graph, $\Delta(G)$ denotes the maximum degree of $G$,
 and if $X\subseteq V(G)$, we denote by $G|X$ 
the subgraph of $G$ induced on $X$.
If $X=\emptyset$, then $\Delta(G|X)=0$. )

\begin{thm}\label{mainthm}
For all integers $t\ge 0$ there is an integer $s$, such that for every graph $G$, if $K_{t+1}$ is not a minor of $G$,
then $V(G)$ can be partitioned into $t$ sets $X_1\l X_t$, such that $\Delta(G|X_i)\le s$ for $1\le i\le t$.
\end{thm}

Such partitions (into subgraphs with bounded maximum degree) are called ``defective colourings'' in the literature --
see for instance~\cite{cowen}.

A reason for interest in \ref{mainthm} is that, despite being much weaker than the original conjecture of Hadwiger, it is still
best possible in the same sense; if we ask for a partition into $t-1$ subgraphs each with bounded maximum degree, the result becomes false.
Let us first see the latter assertion:
\begin{thm}\label{sharp}
For all integers $s\ge0$ and $t\ge 1$, there is a graph $G=G(s,t)$, such that $K_{t+1}$ is not a minor of $G$, and there is no partition
$X_1\l X_{t-1}$ of $V(G)$ into $t-1$ sets such that $\Delta(G|X_i)\le s$ for $1\le i\le t-1$.
\end{thm}
\Proof
If $t=1$ we may take $G(s,t)$ to be a one-vertex graph. For $t\ge 2$, we proceed by induction on $t$. Take the disjoint
union of $s+1$ copies $H_1\l H_{s+1}$ of $G(s,t-1)$, and add one new vertex $v$ adjacent to every other vertex, forming $G=G(s,t)$.
It follows that $G$ has no $K_{t+1}$ minor, since each $H_i$ has no $K_t$ minor.
Assume that $X_1\l X_{t-1}$ is a partition of $V(G)$ into $t-1$ sets such that $\Delta(G|X_i)\le s$ for $1\le i\le t-1$.
We may assume that $v\in X_{t-1}$. If $X_{t-1}\cap V(H_i)\ne \emptyset$ for all $i\in \{1\l s+1\}$, then the degree of $v$ is greater than $s$
in $G|X_{t-1}$, a contradiction; so we may assume that $X_{t-1}\cap V(H_1)=\emptyset$ say. 
Let $Y_i=X_i\cap V(H_1)$ for $1\le i\le t-2$. Then $Y_1\l Y_{t-2}$
provides a partition of $V(H_1)$ into $t-2$ sets; and since $H_1$ is isomorphic to $G(s,t-1)$, it follows that 
$\Delta(H_1|Y_i)> s$ for some $i\in \{1\l t-2\}$, a contradiction. Thus there is no such partition $X_1\l X_{t-1}$. This proves
\ref{sharp}.~\bbox

A warning: it is tempting to view \ref{mainthm} as supporting evidence for Hadwiger's conjecture. However, it is to the same degree ``supporting evidence''
for the false conjecture of Haj\'os~\cite{catlin}, that every graph that contains no subdivision of $K_{t+1}$ is $t$-colourable; because we could
replace the hypothesis of \ref{mainthm} that $G$ has no $K_{t+1}$ minor by the weaker hypothesis that no subgraph of $G$ is a
subdivision of $K_{t+1}$, and the same proof (using an appropriate modification of \ref{extremal}) still works.

\section{The proof}

To prove \ref{mainthm} we use the following lemma, due to Kostochka~\cite{Kostochka1, Kostochka2} and 
Thomason~\cite{Thomason1, Thomason2}. 
\begin{thm}\label{extremal}
There exists $C>0$ such that for all integers $t\ge 0$ and all graphs $G$, if $K_{t+1}$ is not a minor of $G$ then
$G$ has at most $C(t+1)(\log (t+1))^\frac12|V(G)|$ edges.
\end{thm}
We use that to prove two more lemmas:
\begin{thm}\label{smallind}
  Let $t\ge 0$ be an integer, let $C$ be as in \ref{extremal}, and let $r\ge C(t+1)(\log (t+1))^\frac12$.
  Let $G$ be a graph such that $K_{t+1}$ is not a minor of $G$, and
  let $A\subseteq V(G)$ be a stable set of vertices each of degree at least $t$. Then
\[
|E(G\setminus A)|+|A|\le r|V(G\setminus A)|.
\]
\end{thm}

\Proof 
We proceed by induction on $|A|$. By \ref{extremal}, we may assume that $A\neq\emptyset$. Let $v\in A$. Since $v$ has degree at least $t$
and $G$ has no $K_{t+1}$ subgraph, $v$ has two neighbours $x,y$ which are non-adjacent to each other. 
Let $G'=(G\setminus v)+xy$ and $A'=A\setminus\{v\}$. Since $G'$ is a minor of $G$ and so $K_{t+1}$ is not a minor of $G'$, it follows
from the inductive hypothesis that 
$|E(G'\setminus A')|+|A'|\le r |V(G'\setminus A')|=r|V(G\setminus A)|$. 
But $|E(G'\setminus A')|=|E(G\setminus A)|+1$ and $|A'|=|A|-1$. 
This proves~\ref{smallind}.~\bbox

\begin{thm}\label{lowdeg}
Let $t\ge 0$ be an integer, let $C$ be as in \ref{extremal}, and let $r\ge C(t+1)(\log (t+1))^\frac12$ and $r>t/2$.
Let $s$ be the least integer greater than $r(2r-t+2)$.
Let $G$ be a non-null graph, such that $K_{t+1}$ is not a minor of $G$. Then either
\begin{itemize}
\item some vertex has degree less than $t$, or
\item there are two adjacent vertices, both with degree less than $s$.
\end{itemize}
\end{thm}
\Proof We may assume that $t\ge 2$, for if $t\le 1$ the result is trivially true.
Let $A$ be the set of all vertices with degree less than $s$, and $B=V(G)\setminus A$. 
We may assume that every vertex in $A$ has degree at least $t$,
for otherwise the first outcome holds. 
We may also assume that no two vertices of $A$ are adjacent because otherwise the second outcome holds. 
Consequently, by summing all the degrees, we deduce that
$2|E(G)|\ge t|A|+s|B|$.
On the other hand, by \ref{extremal}, $|E(G)|\le r(|A|+|B|)$. It follows that
$t|A|+s|B|\le 2r(|A|+|B|)$, that is, 
$$|A|\ge \frac{s-2r}{2r-t}|B|,$$
since $2r>t$.
But by \ref{smallind}, $|A|\le r|B|$.
Since $G$ is a non-null graph, $|B|\neq 0$ and so $r\ge (s-2r)/(2r-t)$, that is, $s\le r(2r-t+2)$, a contradiction. This proves \ref{lowdeg}.~\bbox

Now we prove \ref{mainthm}, in the following sharpened form.
\begin{thm}\label{mainthm2}
Let $t\ge 0$ be an integer, and let $s$ be as in \ref{lowdeg}. For every graph $G$, if $K_{t+1}$ is not a minor of $G$,
then $V(G)$ can be partitioned into $t$ sets $X_1\l X_t$, such that $\Delta(G|X_i)<s$ for $1\le i\le t$.
\end{thm}
\Proof
We proceed by induction on $|V(G)|+|E(G)|$. 
If some vertex $v$ of $G$ has degree less than $t$, the result follows from the inductive hypothesis by deleting $v$ (find
a partition by induction and add $v$ to some set $X_i$ that contains no neighbour of $v$). If some edge $e$ has both ends
of degree at most $s$, then the result follows from the inductive hypothesis by deleting $e$ (find a partition by induction,
and note that replacing $e$ will not cause either of the ends of $e$ to have degree too large). Thus the result 
follows from \ref{lowdeg}. This proves \ref{mainthm2} and hence \ref{mainthm}.~\bbox

\section{Remark}

Kawarabayashi and Mohar~\cite{kawa} proved the following.
\begin{thm}\label{KM}
There is a function $f(t)$ such that, if $G$ is a graph with no $K_{t+1}$ minor, then $V(G)$ can be partitioned 
into $f(t)$ sets, inducing subgraphs in which every component is of bounded size. 
\end{thm}
Kawarabayashi and Mohar proved that taking $f(t) = \lceil 15.5(t+1)\rceil$ works; and 
Wood~\cite{wood} improved this, showing that taking $f(t)=\lceil 3.5t+2\rceil$ works, using
an unpublished result of Norin and Thomas on large $(t+1)$-connected
graphs with no $K_t$ minor (announced about 2008).
(This has recently been improved to $f(t)=2(t+1)$ by Norin [unpublished].)
That suggests a nice open question -- can we prove the same with $f(t)=t$? This would then give a common extension of
these results and \ref{mainthm}.

Here is a way to improve the Kawarabayashi-Mohar result, showing that
$f(t)=4t$ works (not quite as good as Wood's result, but 
easier).
Alon et al.~\cite[Theorem 6.6]{alon} proved that for all integers $t\ge 0$ and $\Delta$, there exists $s$ such that for every 
graph $G$, if $K_{t+1}$ is not a minor of $G$ and $\Delta(G)\le \Delta$, then $V(G)$ can be partitioned into four sets 
$X_1,\ldots,X_4$ such that every component  of $G|X_1$, $\ldots$, $G|X_4$ has at most $s$ vertices. By combining this 
with \ref{mainthm}, we obtain  a partition of $V(G)$ into $4t$ sets each inducing a graph with no large component.


\begin{thebibliography}{10}

\bibitem{alon}
N.~Alon, G.~Ding, B.~Oporowski, and D.~Vertigan,
\newblock ``Partitioning into graphs with only small components'',
\newblock {\em J. Combin. Theory, Ser. B}, 87(2) (2003), 231--243.

\bibitem{catlin}
P.~A. Catlin,
\newblock ``Haj{\'o}s' graph-coloring conjecture: variations and counterexamples'',
\newblock {\em J. Combin. Theory, Ser. B}, 26(2) (1979), 268--274.

\bibitem{cowen}
L.~Cowen, W.~Goddard, and C.~E. Jesurum,
\newblock ``Defective coloring revisited'',
\newblock {\em J. Graph Theory}, 24(3) (1997), 205--219.

\bibitem{hadwiger}
H.~Hadwiger,
\newblock ``\"{U}ber eine {K}lassifikation der {S}treckenkomplexe'',
\newblock {\em Vierteljschr. Naturforsch. Ges. Z\"urich}, 88 (1943), 133--142.

\bibitem{kawa}
K.~Kawarabayashi and B.~Mohar,
\newblock ``A relaxed {H}adwiger's conjecture for list colorings'',
\newblock {\em J. Combin. Theory, Ser. B}, 97(4) (2007), 647--651.

\bibitem{Kostochka1}
A.~V. Kostochka,
\newblock ``The minimum {H}adwiger number for graphs with a given mean degree of
  vertices'',
\newblock {\em Metody Diskret. Analiz.}, 38 (1982), 37--58.

\bibitem{Kostochka2}
A.~V. Kostochka,
\newblock ``Lower bound of the {H}adwiger number of graphs by their average
  degree'',
\newblock {\em Combinatorica}, 4(4) (1984), 307--316.

\bibitem{RST}
N.~Robertson, P.~Seymour, and R.~Thomas,
\newblock ``Hadwiger's conjecture for {$K_6$}-free graphs'',
\newblock {\em Combinatorica}, 13(3) (1993), 279--361.

\bibitem{Thomason1}
A.~Thomason,
\newblock ``An extremal function for contractions of graphs'',
\newblock {\em Math. Proc. Cambridge Philos. Soc.}, 95(2) (1984), 261--265.

\bibitem{Thomason2}
A.~Thomason,
\newblock ``The extremal function for complete minors'',
\newblock {\em J. Combin. Theory, Ser. B}, 81(2) (2001), 318--338.

\bibitem{wood}
D.~R. Wood,
\newblock ``Contractibility and the {H}adwiger conjecture'',
\newblock {\em European J. Combin.}, 31(8) (2010), 2102--2109.

\end{thebibliography}
\end{document}